\documentclass[twoside, 11pt]{article}
\usepackage{mathrsfs}

\usepackage{mathrsfs,amsfonts,amsmath}
\RequirePackage{ifthen,calc}
\usepackage{theorem}
\usepackage[dvips]{color}

 \catcode`@=11
 \def\@evenhead{\hbox to\textwidth{\footnotesize\rm\thepage \hfill
  {\it }}} 

 \def\@oddhead{\hbox to \textwidth{\footnotesize{\it
  Branching random walk in random environment with random absorption wall  } \hfill\thepage}}

 \renewcommand{\section}{\makeatletter
 \renewcommand{\@seccntformat}[1]{{\csname the##1\endcsname.}\hspace{0.45em}}
 \makeatother \@startsection
{section}
{1}
{0pt}
{\baselineskip}
{0.5\baselineskip}
{\normalsize\bfseries\mathversion{bold}}}

\catcode`@=12
\newcommand\ack{\section*{Acknowledgement}}

\newtheorem{thm}{\noindent Theorem}[section]

\newtheorem{cor}{\noindent Corollary}[section]

\theoremheaderfont{\normalfont\bfseries} \theorembodyfont{\slshape}
{\theorembodyfont{\rmfamily}}
{\theorembodyfont{\rmfamily}}
{\theorembodyfont{\rmfamily}}
{\theorembodyfont{\rmfamily}}
{\theorembodyfont{\rmfamily}}
{\theorembodyfont{\rmfamily}\newtheorem{rem}{\noindent Remark}[section]}
{\theorembodyfont{\rmfamily}}
{\theorembodyfont{\rmfamily}}
{\theorembodyfont{\rmfamily}}

 \def\beqlb{\begin{eqnarray}}\def\eeqlb{\end{eqnarray}}
 \def\beqnn{\begin{eqnarray*}}\def\eeqnn{\end{eqnarray*}}
 \newcommand{\bgeqn}{\begin{equation}}
\newcommand{\edeqn}{\end{equation}}
\def\ra{\rightarrow}

 \numberwithin{equation}{section}
\def\qed{\hfill$\square$\smallskip}
\def\no{\nonumber}

\def\ra{\rightarrow}

\def\bfE{{\mathbb{E}}}

\def\mbfE{{\mathbf{E}}}

\def\bfP{{\mathbb{P}}}

\def\mbfP{{\mathbf{P}}}
\def\bfR{{\mathbb{R}}}

\def\bfN{{\mathbb{N}}}
\def\1{{\mathbf{1}}}

\begin{document}
\title{\bf Quenched small deviation for the trajectory of a random walk with time-inhomogeneous random environment
}
\author{You Lv\thanks{Email: lvyou@dhu.edu.cn },~~Wenming Hong\thanks{Corresponding author. Email: wmhong@bnu.edu.cn}
\\
\\ \small College of Science, Donghua University,
\\ \small Shanghai 201620, P. R. China.
\\
\\ \small School of Mathematical Sciences $\&$ Laboratory of Mathematics and Complex Systems,
\\ \small Beijing Normal University, Beijing 100875, China
}
\date{}
\maketitle


\noindent\textbf{Abstract}: We consider the small deviation probability for random walk with time-inhomogeneous random environment. Compared with the result in Mogul'ski\u{\i} (1974) for the i.i.d. random walk, the rate is smaller (due to the random environment), which is specified in terms of the quenched  and annealed variance.

\smallskip


\smallskip

\smallskip



\section{Introduction for the small deviation principle}

The small deviation problem focuses on the probability that a stochastic process  has fluctuations below its natural scale. We refer to \cite{Mog1974} as one of the well-known results for the classical case. Let $\{X_n\}_{n\in\bfN^+}$ be i.i.d. centered random variables with finite variance $\sigma^2$ and $\{S_n\}_{n\in\bfN}$ the partial sums corresponding to $\{X_n\}_{n\in\bfN^+}$. In a simple but typical case, Mogul'ski\u{\i} \cite{Mog1974} showed that for any $a<0<b,~\alpha\in(0,\frac{1}{2}), x\in\bfR,$
\beqlb\label{1.1Mo}\lim\limits_{n\rightarrow +\infty}\frac{\log\bfP(\forall_{i\leq n}S_i\in[an^{\alpha},bn^{\alpha}]|S_0=x)}{n^{1-2\alpha}}=\frac{-\pi^{2}\sigma^2}{2(b-a)^2}.\eeqlb
Since the random walk has fluctuations of the scale $\sqrt{n},$ the setting $\alpha\in(0,\frac{1}{2})$ means that \eqref{1.1Mo} is of the small deviation principle.

For many years, the small deviation problem has attracted constant attention for a variety of applications in probability such as the Chung type and Strassen type law of iterated logarithm (see \cite{D1983} and \cite{T1993}), fractal structure of fractional Brownian motions (see \cite{X1998}), and the barrier problem of branching random walk (see \cite{GHS2011} and \cite{BJ2012}).  There are also some connections to other topics like the metric entropy of the operator related to a stochastic process (see \cite{KL1993} and \cite{LL1999}). The small deviation problem has been well-studied for many classes of processes like sums of independent random variables (see \cite{Mog1974}, \cite{BM1991} and \cite{QM1995}), Gaussian process and stable processes (see \cite{L1996}, \cite{LQ2004} and the references mentioned therein), fractional Brownian motion and correlated stationary centered Gaussian sequences (see \cite{LL1998} and  \cite{AL2017}), some iterated processes including compound renewal process (see \cite{MFS2008} and \cite{F2013}). 
In the present paper, we want to extend \eqref{1.1Mo} to the random walk with time-inhomogeneous random environment.

Let us start with the definition of this model. We denote $\mu:=\{\mu_n\}_{n\in\bfN^+}$ an i.i.d. sequence with values in the space of probability measures on $\bfR.$ Conditioned to a realization 
of $\mu,$ we sample $\{X_n\}_{n\in\bfN^+}$ a sequence of independent random variables such that for every $n\in\bfN^+,$ the law of $X_n$ is the realization of $\mu_n.$
Set \beqlb\label{sec-sn} S_0=x\in\bfR,~~S_n:=S_0+\sum_{i=1}^{n}X_i.\eeqlb We call $\{S_n\}_{n\in\bfN}$ the {\it random walk with time-inhomogeneous random environment $\mu$ in time}, which is often abbreviated as RWre in the rest of this paper. 
There are two laws to be considered for the RWre. 
We write $\mbfP_{\mu}$ for the law of  $\{S_n\}_{n\in\bfN}$ conditionally on the sequence $\{\mu_n\}_{n\in\bfN^+},$ which is usually called a quenched law. Let $\mathbf{P}$ be the joint law of $\{S_n\}_{n\in\bfN}$ and the environment $\mu,$ which is usually called an annealed law.
The corresponding expectations are denoted by $\mbfE_{\mu}$ and $\mathbf{E}$ respectively. In addition, in order to distinguish from the notation $\mathbf{P}, \mathbf{P}_{\mu}$ in the case of random environment, we should declare that throughout the present paper, we write $\bfP, \bfE$ to present the probability and the corresponding expectation when we do not involve the random environment, just like the notation in \eqref{1.1Mo} and Theorem \uppercase\expandafter{\romannumeral1}. 

Note that the process we consider in the present paper is not the classical
 random walk in random environment which has been well-studied in Zeitouni \cite{Z2004} and many other papers. 
 For the classical random walks in random environment, the random environment is either purely spatial or space-time. However, in our model, the random environment is attached to the time axis. Hence in the annealed sense, the RWre in the present paper can be viewed as the sums of i.i.d. random variables if the random environment $\{\mu_n\}_{n\in\bfN^+}$ is i.i.d.; while in the quenched sense, it is the sums of independent, not necessarily identically distributed random variables. Obviously, if the random environment is degenerate, then $\{S_n\}_{n\in\bfN}$ will degenerate to an i.i.d. centered random walk.



Mallein and Mi{\l}o\'{s} \cite{MM2015} have considered the RWre.
They studied the asymptotics of the probability for the process to stay above some given boundary during $n$ units of time. In a simple but typical case, under the assumptions:
\beqlb\label{asbastien}\exists c',c''>0, \mbfE(e^{c'|M_1|})<+\infty, \mbfE_{\mu}(e^{c'|U_1|})\leq c'',\rm \mathbf{P}-a.s.\eeqlb and (H1),
(where (H1) and the definitions of $M_1, U_1$ are described in the next section,)
 ~they proved the random ballot theorem that \beqlb\label{1.2RBT}\lim\limits_{n\ra \infty}\frac{\log\mbfP_{\mu}(\forall_{i\leq n} S_i\geq 0|S_0=x)}{\log n}=-\gamma_h,~~ \forall x>0, ~~\rm \mathbf{P}-a.s.,\eeqlb where the constant $\gamma_h>1/2$ (resp., $\gamma_h=1/2$) as long as $\mbfP(\mbfE_{\mu}(S_1)=0)<1$ (resp., $\mbfP(\mbfE_{\mu}(S_1)=0)=1$). Note that $\mbfP(\mbfE_{\mu}(S_1)=0)=1$ includes the case that the random environment is degenerate and the well-known conclusion that the  persistence exponent $\gamma_h=1/2$ for an i.i.d. centered random walk with finite variance, from which we can see the impact on the probability
$\mbfP_{\mu}(\forall_{i\leq n} S_i\geq 0|S_0=x)$ from the random environment.

 It should be mentioned that the time-inhomogeneous version of many-to-one formula introduced by Mallein \cite{M2015a} constructs the relationship between RWre and the branching random walk with time-inhomogeneous random environment(BRWre). Based on this relationship, the random ballot theorem is a key tool to study the maximal displacement of a BRWre (see Mallein and Mi{\l}o\'{s} \cite{MM2016}).
The initial motivation to analyzing the problems in the present paper is also from the studying of BRWre. The main result in the present paper can be used to study the barrier problem of a BRWre, which is given in the companion paper \cite{LY3}.


 The rest of the article is organized as follows. Section 2 contains the assumptions and results. Since Theorem 2.2 implies Theorem 2.1 and we can obtain Theorem 2.3 by checking the proof of Theorem 2.2, we only need to give the proofs of Theorem 2.2 and Corollary 2.1 in Section 3. The proof of \eqref{bmt+++} is given in Appendix.

\section{Assumptions and results}

The RWre $\{S_n\}_{n\in\bfN}$ (with the random environment $\mu$) is defined as \eqref{sec-sn}. Denote $$M_{n}:=\mbfE_\mu(S_n), ~U_{n}:=S_n-\mbfE_\mu(S_n),~\Gamma_n:=\mbfE_\mu(U^2_n)=\mbfE_\mu(S^{2}_n)-M^2_{n}.$$ 
Obviously, it is true that $\mbfE(U_n^2)=\mbfE(\Gamma_n).$ We assume that
\begin{itemize}
  \item[(H1)] $\mbfE M_1=0,\sigma^2_{_A}:=\mbfE(M_1^2)\in[0,+\infty),\sigma^2_{_Q}:=\mbfE(U_1^2)\in(0,+\infty).$
  \item[(H2)] There exist $\lambda_1>\lambda_0> 2$ such that $\mbfE(|M_{1}|^{\lambda_1})<+\infty.$
  \item[(H3)] There exist $\lambda_2> 2, \lambda_3>\max\{\frac{\lambda_0}{2}, 2\}$ such that $\mbfE\left(\left[\mbfE_\mu(|U_1|^{\lambda_2})\right]^{\lambda_3}\right)<+\infty.$
\end{itemize}
\begin{rem}
The $\lambda_0$ is set for a more concise expression in the proof and its value will be given in Theorem 2.1-2.2 and Corollary 2.1. We should note that the role of $\lambda_2$ is different from that of $\lambda_1$ and $\lambda_3$ in some extent. $\lambda_2$ represents the  integrability in quenched sense, while $\lambda_1$ and $\lambda_3$ represent the integrability in annealed sense. In other words, we can see that (H2) and (H3) will  automatically hold for any large $\lambda_1$ and $\lambda_3$ when $\mu_1$ is degenerate (or even has only finite states). Moreover, we can see that the assumptions (H2) and (H3) are weaker than the assumption \eqref{asbastien} given in \cite{MM2016}.

\end{rem}

In order to obtain an union bound, 
we need to consider the shift of the starting time since $\{S_n\}_{n\in\bfN}$ is a time-inhomogeneous random walk under $\mbfP_{\mu}$.
\begin{thm}\label{mogbefore}
Let $a<0<b, x\in\bfR.$ Under the assumptions (H1)-(H3), for any $\alpha\in(0,\frac{1}{2})$,
the convergence
\beqlb\label{upper}\lim\limits_{n\rightarrow +\infty}\frac{\log \mbfP_\mu
(\forall_{k\leq i\leq k+n} S_{i}\in[an^\alpha, bn^\alpha]|S_{k}=x)}{n^{1-2\alpha}}= \frac{-\sigma^2_Q}{(b-a)^2}\gamma\left(\frac{\sigma_A}{\sigma_Q}\right)\eeqlb
holds in the sense of ${\rm \mathbf{P}-a.s.}$ if $\lambda_0=2\alpha^{-1}$ and in Probability $\mathbf{P}$ if $\lambda_0=\alpha^{-1},$ where the $\lambda_0$ has been described in Assumptions (H2)-(H3) and the real-valued function $\gamma$ is defined in the Theorem \uppercase\expandafter{\romannumeral1} below.\end{thm}


\noindent\emph{{\bf Theorem \uppercase\expandafter{\romannumeral1}~(\cite[Theorem 2.1 and Theorem 3.1]{LY201801})} Let $B,W$ be two independent standard Brownian motions, $W_0=0, B_0=x$. Let the real constants $a,b,a_0,b_0,a',b',\beta$ satisfy $a<a_0\leq b_0<b$ and $a\leq a'< b'\leq b.$ 
Define
$\underline{\mathcal{X}}_t:=-\log\sup\limits_{x\in\bfR} \bfP(\forall_{s\in [0,t]} B_{s}+\beta W_{s}\in[a,b]|W);$
$\overline{\mathcal{X}}_t:=-\log\inf\limits_{x\in[a_0,b_0]} \bfP(\forall_{s\in [0,t]} B_{s}+\beta W_{s}\in[a,b],~B_{t}+\beta W_{t}\in[a',b']|W).$
Then we have\begin{eqnarray}\label{bmt+}
\lim\limits_{t\rightarrow+\infty}\frac{\underline{\mathcal{X}}_t}{t}=\lim\limits_{t\rightarrow+\infty}\frac{\overline{\mathcal{X}}_t}{t}=\frac{\gamma(\beta)}{(b-a)^2},~~~{\rm a.s. ~~~and~~ in}~ L_p~(p\geq1),
\end{eqnarray}
where the $\gamma$ is a well-defined, positive, convex and even function with
\beqlb\label{bmt1}\gamma(0)=\frac{\pi^2}{2},~~ \gamma(\beta)\geq\frac{\pi^2(1+\beta^2)}{2}, ~~\forall \beta\in\bfR.\eeqlb Moreover, for any $t>0, d>0$ and $\tilde{d}>1,$ it is true that
\begin{eqnarray}\label{bmt++}\lim\limits_{n\rightarrow +\infty}e^{dn} \bfP\big(\overline{\mathcal{X}}_t \geq  (\log n)^{\tilde{d}}\big)=0.\end{eqnarray}}
\noindent In fact, by checking the proof for \eqref{bmt++} in \cite[Theorem 3.1]{LY201801} carefully, \eqref{bmt++} can be strengthened to that for any given
$\overline{\mathcal{X}}_t,$ 
\begin{eqnarray}\label{bmt+++}\forall d\in\bfR^+,~~~\bfE(e^{d\overline{\mathcal{X}}_t})<+\infty.\end{eqnarray} The proof of \eqref{bmt+++} will be given in Appendix.

\begin{rem}\label{rem2}
(1)~Note that if $\sigma_{_A}=0~({\rm i.e.}~\mbfP(\mbfE_{\mu}(S_1)=0)=1)$, there contains two cases.

\emph{Case 1}: The random environment is degenerate. That is to say, our model will degenerate to an i.i.d. centered random walk and thus $\sigma^2_{_Q}=\mbfE(S^2_1)$. Note that $\gamma(0)=\frac{\pi^2}{2},$ hence in this case \eqref{upper} is consistent with the Mogul'ski\u{\i} Theorem \eqref{1.1Mo}.~On the other hand, \eqref{bmt1} implies $\frac{-\sigma^2_Q}{(b-a)^2}\gamma\left(\frac{\sigma_A}{\sigma_Q}\right)\leq\frac{-\pi^2\mathbf (\sigma^2_A+\sigma^2_Q)}{2(b-a)^2}.$ Comparing this fact with \eqref{1.1Mo}, we can see that the quantity $\sigma^2_A$ reflects the impact on the asymptotic behavior of the small deviation probability from the random environment.

\emph{Case 2}: The random environment is non-degenerate and $\sigma_{_A}=0.$ That is to say, conditionally on any realization of $\mu,$ $\{S_n\}_{n\in\bfN}$ is an independent, not necessarily identically distributed (i.n.d.) centered random walk. For $c\in\bfR^+$, by our result one can see
$$\lim\limits_{n\rightarrow +\infty}\frac{\log \mbfP_\mu
(\forall_{k\leq i\leq k+n} S_{i}\in[-cn^\alpha, cn^\alpha]|S_k=0)}{n^{1-2\alpha}}= \frac{-\sigma^2_{_Q}}{4c^2}\gamma\left(0\right)=\frac{-\pi^2\sigma^2_{_Q}}{8c^2},~~~{\rm \mathbf{P}-a.s.}.$$
In addition, according to Remark 2.1(2) and the Borel Cantelli 0-1 law,  we can obtain that $\lim\limits_{n\rightarrow +\infty}n^{-1}\mbfE_{\mu}(|U_{k+n}-U_k|^2)=\sigma^2_{_Q}, \forall k\in\bfN,~{\rm \mathbf{P}-a.s.}.$ Note that $S_n=U_n$ when $\sigma_A=0,$ hence we
can rewrite the above convergence as
\beqlb\label{add37}\lim\limits_{n\rightarrow +\infty}\frac{(cn^{\alpha})^2\log \mbfP_\mu
(\forall_{k\leq i\leq k+n} S_{i}\in[-cn^\alpha, cn^\alpha]|S_k=0)}{n\cdot n^{-1}\mbfE_{\mu}(|S_{k+n}-S_k|^2)}=\frac{-\pi^2}{8},~~~{\rm \mathbf{P}-a.s.}.\eeqlb

On the other hand, let us recall the well-known result on the small deviation probability for the sums of i.n.d. random variables in Shao \cite{QM1995}. For a random sequence $\{T_n\}_{n\in\bfN}$ with independent increment satisfying $\bfE(T_n)=0, \bfE(T^2_n)<+\infty, \forall n\in\bfN$ and an uniform Lindeberg's condition, Shao showed that
\beqlb\label{shao}\lim\limits_{n\rightarrow +\infty}\frac{x^2_{n,k}\log\bfP(\forall_{i\leq n} |T_{k+i}-T_k|\leq x_{n,k})}{\bfE(|T_{k+n}-T_k|^2)}=-\frac{\pi^2}{8}, ~\forall k\in\bfN.\eeqlb
 as long as $\{x_{n,k}\}_{n\in\bfN, k\leq n}$ satisfies that $\lim\limits_{n\rightarrow +\infty}\min\limits_{k\in\bfR}x_{n,k}=\lim\limits_{n\rightarrow +\infty}\min\limits_{k\in\bfR}x^{-2}_{n,k}\bfE(|T_{k+n}-T_k|^2)=+\infty.$

 Now let us set $x_{n,k}:=cn^{\alpha}, \alpha\in(0,1/2), c>0.$ By the above analysis we can see our result \eqref{add37} is consistent with Shao Theorem \eqref{shao}.



(2)~Of course, for Theorem 2.1, the more noteworthy part is the case that $\sigma_A>0,$ which is the typical case for RWre. That is to say, conditionally on any realization of $\mu,$ $\{S_n\}_{n\in\bfN}$ is an i.n.d. and not necessarily centered random walk. One the whole, for the case $\sigma_A>0,$ Shao's result can show the asymptotics of the probability for a random walk staying in a neighborhood of its expectations while we show the probability for a random walk staying in a neighborhood of $0.$ Let us show the details as follows.

Since conditionally on any realization of $\mu,$ ${S_{i}-\mbfE_{\mu}(S_{i})}$ satisfies the conditions in Shao Theorem $\eqref{shao}$, we can apply $\eqref{shao}$ to RWre to get the conclusion that
\beqlb\label{shao+}\lim\limits_{n\rightarrow +\infty}\frac{(cn^{\alpha})^2\log\mbfP_{\mu}(\forall_{i\leq n} |S_{i}-\mbfE_{\mu}(S_{i})|\leq cn^{\alpha})}{\mbfE_{\mu}(|S_n-\mbfE_{\mu}(S_{n})|^2)}=-\frac{\pi^2}{8},~{\rm \mathbf{P}-a.s.}.\eeqlb
Note that $\lim\limits_{n\rightarrow +\infty}n^{-1}\mbfE_{\mu}(|S_n-\mbfE_{\mu}(S_{n})|^2)=\sigma^2_{_Q},~{\rm \mathbf{P}-a.s.}.$ Hence \eqref{shao+} is equal to say \beqlb\label{shao++}\lim\limits_{n\rightarrow +\infty}\frac{\log\mbfP_{\mu}(\forall_{i\leq n} S_{i}\in [\mbfE_{\mu}(S_{i})-cn^{\alpha},\mbfE_{\mu}(S_{i})+cn^{\alpha}])}{n^{1-2\alpha}}=-\frac{\pi^2\sigma^2_{_Q}}{8c^2},~{\rm \mathbf{P}-a.s.}.\eeqlb
On the other hand, Theorem 2.1 shows that $\lim\limits_{n\rightarrow +\infty}\frac{\log\mbfP_{\mu}(\forall_{i\leq n} S_{i}\in [-cn^{\alpha},cn^{\alpha}])}{n^{1-2\alpha}}=-\frac{\sigma^2_Q}{4c^2}\gamma\left(\frac{\sigma_A}{\sigma_{Q}}\right),~{\rm \mathbf{P}-a.s.}.$ Moreover, by the property of function $\gamma$ we can see that $-\frac{\sigma^2_Q}{4c^2}\gamma\left(\frac{\sigma_A}{\sigma_Q}\right)<-\frac{\pi^2\sigma^2_Q}{8c^2}$ only if $\sigma_A>0.$ The strict inequality reflects that for our model, the probability that the trajectory of $\{S_i\}_{i\leq n}$ stays in the $cn^{\alpha}$-neighborhood of its quenched mean $\{\mbfE_{\mu}(S_i)\}_{i\leq n}$ is larger than that in the $cn^{\alpha}$-neighborhood of its annealed mean $0.$

(3)~Note that (1.1) only need a finite second moment. Therefore, for the random environment case,  a future work which may be worth considering is what will happen if $\lambda_2=2$ in assumption (H3). We believe this critical case may contain more interesting phenomenon.
\end{rem}

While proving Theorem \ref{mogbefore}, we actually obtain a stronger result. We can further show the uniformity of the estimate with respect
to the starting position and ending position. Moreover, the convergence in Theorem \ref{mogbefore} remains true when the starting time depends on $n$.
For the sake of simplicity, throughout this paper, we denote \beqlb\label{bmt2}\tilde{\gamma}:=\sigma^2_Q\gamma\left(\frac{\sigma_A}{\sigma_Q}\right).\eeqlb The refined version can be stated as follows.
\begin{thm}\label{mog}
Let $\{t_n\}_{n\in\bfN}$ be a sequence of non-negative integers 
and $\bar{t}_n:= t_n+n.$~Set $a<a_0\leq b_0<b,a\leq a'<b'\leq b.$ 
 Under the assumptions (H1)-(H3), for any $\alpha\in(0,\frac{1}{2}),$ the following convergence
\beqlb\label{upper0}\lim\limits_{n\rightarrow +\infty}\frac{\sup\limits_{x\in\bfR}\log \mbfP_\mu
\left(\forall_{t_n\leq i\leq \bar{t}_n} S_{i}\in[an^\alpha, bn^\alpha]|S_{t_n}=x\right)}{n^{1-2\alpha}}=\frac{-\tilde{\gamma}}{(b-a)^2},\eeqlb
\beqlb\label{lower0}&&
\lim\limits_{n\rightarrow +\infty}\frac{\inf\limits_{x\in[a_0 n^{\alpha}, b_0 n^{\alpha}]}\log \mbfP_\mu
\left(\begin{split}\forall_{t_n\leq i\leq \bar{t}_n} S_{i}\in[an^\alpha, bn^\alpha],\\~S_{\bar{t}_n}\in [a'n^\alpha,b'n^\alpha]\end{split}\Bigg|S_{t_n}=x\right)}{n^{1-2\alpha}}=\frac{-\tilde{\gamma}}{(b-a)^2}\eeqlb
hold in the sense of ~${\rm\mathbf{P}-a.s.}$ if $\lambda_0=2\alpha^{-1}$ and in Probability $\mathbf{P}$ if $\lambda_0=\alpha^{-1}.$
\end{thm}


In addition, checking the proof of Theorem \ref{mog} in the next section,  we see the convergence will also hold for two independent random walks, which is stated in the following theorem.
\begin{thm}\label{mog++++}
Let $V, \hat{V}$ be two independent centered i.i.d. random walks with $\hat{V}_0=V_0=0,~\bfE(|\hat{V}_1|^{q_1})+\bfE(|V_1|^{q_2})<+\infty,~ q_1>2, \bfE(\hat{V}_1^{2})>0$. Denote $V_{i}^{(k)}=V_{i+k}-V_k, \hat{V}_{i}^{(k)}=\hat{V}_{i+k}-\hat{V}_k$. The setting of $\alpha, a,a_0,a',b,b_0,b',t_n$ are the same as what we describe in Theorem 2.2.
~Denote $\sigma^2_{_V}:=\bfE(V_1^2), \sigma^2_{_{\hat{V}}}:=\bfE(\hat{V}_1^2), \hat{\gamma}=\sigma^2_{_{\hat{V}}}\gamma\left(
\frac{\sigma_{_{V}}}{\sigma_{_{\hat{V}}}}\right).$
Then for any $\alpha\in(0, \frac{1}{2}),$ the following convergence
$$\lim\limits_{n\ra\infty}\frac{\sup\limits_{x\in\bfR}\log\bfP
(\forall_{i\leq n} ~x+\hat{V}^{(t_n)}_i\in[an^\alpha+V^{(t_n)}_i, bn^\alpha+V^{(t_n)}_i]|V)}{n^{1-2\alpha}}=\frac{-\hat{\gamma}}{(b-a)^2},
$$
$$\lim\limits_{n\ra\infty}\frac{\inf\limits_{x\in[a_0n^{\alpha},b_0n^{\alpha}]}\log\bfP
\left( \begin{split}\forall_{i\leq n} ~x+\hat{V}_i^{(t_n)}\in[an^\alpha+V_i^{(t_n)}, bn^\alpha+V_i^{(t_n)}]\\x+\hat{V}_{n}^{(t_n)}\in[a'n^\alpha+V^{(t_n)}_{n}, b'n^\alpha+V^{(t_n)}_{n}]\end{split}\Bigg|V\right)}{n^{1-2\alpha}}=\frac{-\hat{\gamma}}{(b-a)^2}
$$
hold in the sense of ${\rm \bfP-a.s.}$ if $q_2>2\alpha^{-1}$ and in Probability $\bfP$ if $q_2>\alpha^{-1}.$
\end{thm}

In fact, Mogul'ski\u{\i} \cite{Mog1974} considered the small deviation probability with a more general boundary. Therefore, we present that the $a, b$ in Theorem 2.2 can be extended to two continuous functions, which is given as the following corollary.
\begin{cor}\label{mogc1}
Let $g(s), h(s)$ be two continuous functions on $[0,1]$ and $g(s)<h(s)$ for any $s\in [0,1].$ We set $g(0)< a_0\leq b_0 <h(0), g(1)\leq a'<b'\leq h(1)$
~and~ $C_{g,h}:=\int_{0}^{1}\frac{1}{[h(s)-g(s)]^2}ds.$ The setting of $t_n, \bar{t}_n$ are the same as what we describe in Theorem 2.2. Then for any $\alpha\in (0,\frac{1}{2}),$ under the assumptions (H1)-(H3), the convergence
\beqlb\label{mogc1u}\lim\limits_{n\rightarrow +\infty}\frac{\sup\limits_{x\in\bfR}\log\mbfP_\mu
\Big(\forall_{t_n\leq i\leq \bar{t}_n}S_{i}\in \big[g\big(\frac{i-t_n}{n}\big)n^\alpha,h\big(\frac{i-t_n}{n}\big)n^\alpha\big]\Big|S_{t_n}=x\Big)}{n^{1-2\alpha}}= -C_{g,h}\tilde{\gamma},\eeqlb
\beqlb\label{mogc1l}&&\lim\limits_{n\rightarrow +\infty}\frac{\inf\limits_{x\in[a_0 n^{\alpha}, b_0 n^{\alpha}]}\log \mbfP_\mu
\left(\forall_{t_n\leq i\leq \bar{t}_n} S_{i}\in \big[g\big(\frac{i-t_n}{n}\big)n^\alpha,h\big(\frac{i-t_n}{n}\big)n^\alpha\big], S_{\bar{t}_n}\in [a'n^\alpha,b'n^\alpha] \Big|S_{t_n}=x\right)}{n^{1-2\alpha}}\no
\\&&~~~~~~~~~~~~~~~~~~~~~~~~~~~~~~~~~~~~~~~~~~= -C_{g,h}\tilde{\gamma}\eeqlb
hold in the sense of ${\rm \mathbf{P}-a.s.}$ if $\lambda_0=2\alpha^{-1}$ and in Probability $\mathbf{P}$ if $\lambda_0=\alpha^{-1}.$
\end{cor}


\section{Proofs}
\subsection{Sakhanenko's theorem}
~~The main idea of the proof is coupling the RWre with two Brownian motions by the main theorems in 
Sakhanenko \cite{Sak2006}, which is a celebrated theorem in strong approximation theory. 
Here we state the theorem below in order to facilitate the applying and understanding on the following proofs.

\noindent\emph{{\bf Theorem \uppercase\expandafter{\romannumeral2}~(Sak, \cite[Theorem 1]{Sak2006})}
Let~$\xi_1,\xi_2,\ldots,\xi_n$ be an infinite sequence of independent random variables satisfying $\bfE(\xi_j)=0, \bfE(\xi^2_j)<+\infty$.
Then for any fixed number $\lambda^{*}\geq 2$ there exists a standard Brownian motion $W$ (depend on $\lambda^{*}$, and the distributions of $\xi_1,\xi_2,\ldots,\xi_n$) such that
$$\forall x\geq y>0, \bfP\left(\max_{k\leq n}|S_k-W_{D_k}|\geq C^{*}\lambda^{*}x\right)\leq \left[\left(y^{-\lambda^{*}}\sum_{k=1}^n\bfE(\xi^{\lambda^{*}}_k)\right)\right]^{x/y}+\bfP(\max_{k\leq n}\xi_k>y),$$
where $D_k:=\sum_{i=1}^k\bfE(\xi^{2}_i)$ and $C^{*}$ is an absolutely constant.}

In particular, if we choose $x=y,$ then we get
\begin{eqnarray}\label{sak06}\forall x>0, \bfP\left(\max_{k\leq n}|S_k-W_{D_k}|\geq C^{*}\lambda^{*}x\right)\leq 2x^{-\lambda^{*}}\sum_{k=1}^n\bfE(\xi^{\lambda^{*}}_k).\end{eqnarray}

In the present paper, we draw lessons from the coupling method, which has been used in Mallein and Mi{\l}o\'{s} \cite{MM2015}. We put the outline of the proof in the next subsection. 

\subsection{The lower bound of Theorem 2.2}

~~In the main line of the proof we let $\lambda_0=2\alpha^{-1}$ and devote to show the ${\rm \mbfP-a.s.}$ part, while the case $\lambda_0=\alpha^{-1}$ will be marked in parentheses when necessary.

To prove Theorem 2.2, we only need to show these two inequalities:
\begin{eqnarray}\label{1+}\varlimsup\limits_{n\rightarrow +\infty}\frac{\sup\limits_{x\in\bfR}\log \mbfP_\mu
\left(\forall_{t_n\leq i\leq \bar{t}_n} S_{i}\in[an^\alpha, bn^\alpha]|S_{t_n}=x\right)}{n^{1-2\alpha}}\leq\frac{-\tilde{\gamma}}{(b-a)^2},~{\rm \mathbf{P}-a.s.;}\end{eqnarray}
\begin{eqnarray}\label{2-}\varliminf\limits_{n\rightarrow +\infty}\frac{\inf\limits_{x\in[a_0 n^{\alpha}, b_0 n^{\alpha}]}\log \mbfP_\mu
\left(\begin{split}\forall_{t_n\leq i\leq \bar{t}_n} S_{i}\in[an^\alpha, bn^\alpha],\\~S_{\bar{t}_n}\in [a'n^\alpha,b'n^\alpha]\end{split}\Bigg|S_{t_n}=x\right)}{n^{1-2\alpha}}\geq\frac{-\tilde{\gamma}}{(b-a)^2},~{\rm \mathbf{P}-a.s.}.\end{eqnarray}
We first show the lower bound.

\noindent{\bf Proof of Theorem \ref{mog}: the lower bound \eqref{2-}}

  Let $D\in\bfN^+, ~T:=\lfloor Dn^{2\alpha}\rfloor,~ K:=\left\lfloor \frac{n}{T} \right\rfloor,$ $t_{n,k}:=t_n+kT,$ where $\lfloor x\rfloor$ represents the maximum integer value not exceeding $x.$ Without loss of generality, in the proof of the lower bound we always assume that $a<a_0\leq a'<b'\leq b_0<b.$ Throughout the proof, we choose $a'',b'',a''',b'''$ such that
  \begin{eqnarray}\label{3.3.1}a'<a''<a'''<b'''<b''<b'. \end{eqnarray}
  Recall that~$\bar{t}_n:=t_n+n.$ By Markov property and some standard calculations, for $n$ large enough, it is true that
\begin{eqnarray}\label{3.3.2}
&&\inf\limits_{x\in[a_0 n^{\alpha}, b_0 n^{\alpha}]}\mbfP_\mu
\Big(\forall_{t_n\leq i\leq \bar{t}_n}~ S_{i}\in[an^{\alpha},bn^{\alpha}],~S_{\bar{t}_n}\in[a'n^{\alpha}, b'n^{\alpha}] \Big|S_{t_n}=x\Big)\no
\\&\geq&\prod\limits_{k=0}\limits^{K-1}\inf\limits_{x\in [a_0 n^\alpha,b_0 n^\alpha]} \mbfP_\mu\left(\forall_{i\leq T}~ S_{t_{n,k}+i}\in[a n^\alpha, bn^\alpha], ~S_{t_{n,k+1}}\in[a''n^\alpha,b''n^\alpha] \Bigg|S_{t_{n,k}}=x\right)\no
\\&\times&\inf\limits_{x\in [a'' n^\alpha,b'' n^\alpha]} \mbfP_\mu\left(\forall_{i\leq \bar{t}_n-t_{n,K}}~ S_{t_{n,K}+i}\in[an^\alpha, b n^\alpha], S_{\bar{t}_n}\in[a'n^\alpha,b'n^\alpha] \Bigg|S_{t_{n,K}}=x\right).
\end{eqnarray}

Let us give a quick sketch of the proof first.  We decompose $S_{n}$ into the sum of $M_n$ and $U_n.$ (Recall that $M_n:=\mbfE_\mu(S_n), ~U_{n}:=S_n-\mbfE_\mu(S_n).$) By Theorem \uppercase\expandafter{\romannumeral2} , we use two independent standard Brownian motions $B$ and $\hat{W}$ to approximate $\{U_n\}_{n\in\bfN}$ and $\{M_n\}_{n\in\bfN}$ respectively such that for large enough $n,$
$$\mbfP_\mu\left(\forall_{i\leq T}S_{t_{n,k}+i}\in[an^\alpha, bn^\alpha]\Big|S_{t_{n,k}}=0\right)\approx\mbfP(\forall_{s\leq D}B_{s\sigma^2_{_Q}}+n^{-\alpha}\hat{W}_{s\sigma_{_A}^2n^{2\alpha}}^{(n,k)}\in[a,b]|\hat{W}^{(n,k)}), \forall k\in\bfN,$$
holds with overwhelming probability (with respect to $\mbfP$), where $\hat{W}^{(n,k)}$ is an independent copy of $\hat{W}$, $a<0<b, B_0=\hat{W}^{(n,k)}_0=0.$ 
Finally we complete the proof by Theorem \uppercase\expandafter{\romannumeral1} and the law of large numbers.

Next, we will divide the proof into four steps.

\emph{step 1: Using a Brownian motion to approximate $\{U_{t_{n,k}+i}-U_{t_{n,k}}\}_{i\leq T}$.}

Recalling that $\Gamma_n:=\mbfE_\mu(U^2_n).$ Let~$U^{(n,k)}_{i}:=U_{t_{n,k}+i}-U_{t_{n,k}}, M^{(n,k)}_i:=M_{t_{n,k}+i}-M_{t_{n,k}}, \Gamma^{(n,k)}_i:=\Gamma_{t_{n,k}+i}-\Gamma_{t_{n,k}}, \zeta_i:=\mbfE_{\mu}((U_{i}-U_{i-1})^{\lambda_2})$ and $\zeta_{n,k}:=\sum_{i=t_{n,k}+1}^{t_{n,k}+T}\mbfE_{\mu}((U_{i}-U_{i-1})^{\lambda_2}).$
 According to Theorem \uppercase\expandafter{\romannumeral2}, for any given $\mu,$ we can construct a standard Brownian motion $B^{(n,k)}$ such that 
\begin{eqnarray}\label{3.3.8}\forall n\in\bfN, ~k\leq K,~\lambda^*\geq 2,~\mbfP_{\mu}\left(\Delta_{n,k}\geq x\right)\leq \frac{2(C^*\lambda^*)^{\lambda^*}\zeta_{n,k}}{x^{\lambda^*}},~~{\rm \mathbf{P}-a.s.,}\end{eqnarray} where $\Delta_{n,k}:=\max_{j\leq T}\left|U^{(n,k)}_{j}-B^{(n,k)}_{\Gamma^{(n,k)}_j}\right|.$
For any $n,k,$ in the following proof we always write $B^{(n,k)}$ as $B$ for simplicity since our deduction only involves the distribution of $B^{(n,k)}$ but not its sample. Choose $\varsigma\in(\max\{\frac{\alpha}{2},\frac{2}{\lambda_1},\frac{2\alpha}{\lambda_2},\frac{1}{\lambda_3}\},\alpha)$ and $\rho\in(\frac{\alpha+\varsigma}{2},\alpha)$ (for the case $\lambda_0=\alpha^{-1},$ we choose $\varsigma\in(\max\{\frac{\alpha}{2},\frac{1}{\lambda_1},\frac{2\alpha}{\lambda_2},\frac{1}{2\lambda_3}\},\alpha)$) and define $\mathcal{V}_{n,k}:=\{\Delta_{n,k}\leq n^{\rho}\}.$ Then we have
\begin{eqnarray}\label{3.3.9}
&&\mbfP_\mu\left(\forall_{i\leq T}~ S_{t_{n,k}+i}\in[a n^\alpha, bn^\alpha], ~S_{t_{n,k+1}}\in[a''n^\alpha,b''n^\alpha] \Big|S_{t_{n,k}}=x\right)\no
\\&=&\mbfP_\mu\left(\forall_{i\leq T} ~ x+U^{(n,k)}_{i}+M^{(n,k)}_i\in[an^\alpha, bn^\alpha], x+U^{(n,k)}_{T}+M^{(n,k)}_T\in[a'' n^\alpha,b'' n^\alpha],\mathcal{V}_{n,k}\right)\nonumber
\\&\geq& \mbfP_\mu\left(\begin{aligned}&\forall_{i\leq T}~ x+B_{\Gamma^{(n,k)}_i}+M^{(n,k)}_i\in[an^\alpha+n^\rho, bn^\alpha-n^\rho],
\\&~~~x+B_{\Gamma^{(n,k)}_T}+M^{(n,k)}_T\in[a'' n^\alpha+n^\rho,b'' n^\alpha-n^\rho]\end{aligned}\right)-\mbfP_\mu(\mathcal{V}^{c}_{n,k}).
\end{eqnarray}
The inequality in \eqref{3.3.9} is due to the definition of the event $\mathcal{V}_{n,k}.$

 Define
\begin{eqnarray}\label{3.3.6}F_n:=\cap_{k=0}^{K}\left\{\forall_{i\leq T}|\Gamma^{(n,k)}_{i}-i\sigma^2_Q|\leq 2C^*\lambda_3n^{\alpha+\varsigma}\right\}, \hat{F}_n:=\cap_{k=0}^{K}\left\{\zeta_{n,k}\leq 2\mbfE(U^{\lambda_2}_1)T\right\}.\end{eqnarray} Recall that $\mbfE(\Gamma^{(n,k)}_i)=i\sigma_Q^2.$
Now on the event $F_n,$ we estimate the difference between $B_{\Gamma^{(n,k)}_i}$ and $B_{s\sigma_Q^2}, s\in[i,i+1).$ Define
\begin{eqnarray}\label{3.3.10}
\Theta_{n}:=\left\{\forall i\leq T,~\forall |t|\leq 2C^*\lambda_3n^{\alpha+\varsigma}+\sigma^2_Q, ~|B_{i\sigma^2_{Q}}-B_{i\sigma^2_{Q}+t}|\leq \frac{1}{2}n^\rho\right\}.
\end{eqnarray}
Note that 
$\left|B_{s\sigma^2_Q}-B_{\Gamma^{(n,k)}_i}\right|\leq \left|B_{s\sigma^2_Q}-B_{i\sigma^2_Q}\right|+\left|B_{i\sigma^2_Q}-B_{\Gamma^{(n,k)}_i}\right|.$
Hence on the event $F_n$, we have
\begin{eqnarray}\label{3.3.11}
&&\mbfP_\mu\left(\begin{aligned} &\forall_{i\leq T}~ x+B_{\Gamma^{(n,k)}_i}+M^{(n,k)}_i\in[an^\alpha+n^\rho,bn^\alpha-n^\rho],
\\&~~ x+B_{\Gamma^{(n,k)}_T}+M^{(n,k)}_T\in[a'' n^\alpha+n^\rho,b'' n^\alpha-n^\rho],~\Theta_n\end{aligned}\right)\nonumber
\\&\geq& \mbfP_{\mu}\left(\begin{aligned}\forall_{i<T,~\lfloor s\rfloor=i} ~x+B_{s\sigma^2_Q}+M^{(n,k)}_i\in[an^\alpha+\frac{3}{2}n^\rho,bn^\alpha-\frac{3}{2}n^\rho],
\\ x+B_{T\sigma^2_Q}+M^{(n,k)}_T\in[a'' n^\alpha+\frac{3}{2}n^\rho,b'' n^\alpha-\frac{3}{2}n^\rho]\end{aligned}\right)-\mbfP_\mu(\Theta^c_{n})\no
\\&=& \mbfP_{\mu}\left(\begin{aligned}\forall_{i<T,~\lfloor s\rfloor=i} ~x+\sigma_Q B_{s}+M^{(n,k)}_i\in[an^\alpha+2n^\rho,bn^\alpha-2n^\rho],
\\ x+\sigma_Q B_{T}+M^{(n,k)}_T\in[a'' n^\alpha+2n^\rho,b'' n^\alpha-2n^\rho]\end{aligned}\right)-\mbfP_\mu(\Theta^c_{n}).
\end{eqnarray}
Thanks to~\eqref{3.3.8}, on the event $\hat{F}_n$ we have
\begin{eqnarray}\label{3.3.12}
\mbfP_\mu(\mathcal{V}^c_{n,k})\leq 4(C^*\lambda_2)^{\lambda_2}\mbfE(U^{\lambda_2}_1)Tn^{-\rho\lambda_2}.
\end{eqnarray}
On the other hand, by the Cs\"{o}rg\H{o} and R\'{e}v\'{e}sz's estimation \cite[Lemma 1]{CR1979} for Brownian motion, we know there exist constants $c_1,c_2>0$ such that
\begin{eqnarray}\label{3.3.13}
\mbfP_\mu(\Theta^c_{n})
\leq2T\cdot\mbfP_{\mu}\left(\sup_{t\in[0,n^{\alpha+\varsigma}+\sigma^2_Q]}|B_t|>n^{\rho}\right)\leq 2T\cdot c_1\exp\left\{\frac{-c_2n^{2\rho}}{n^{\alpha+\varsigma}+\sigma^2_Q}\right\}.
\end{eqnarray}
Combining \eqref{3.3.11}-\eqref{3.3.13} with \eqref{3.3.9}, on the event $F_n,$ we can find a $c_3>0$ such that
\begin{eqnarray}\label{3.3.14}
&&\mbfP_\mu\left(\forall_{i\leq T}~ S_{t_{n,k}+i}\in[a n^\alpha, bn^\alpha], ~S_{t_{n,k+1}}\in[a''n^\alpha,b''n^\alpha] \Big|S_{t_{n,k}}=x\right)\no
\\&\geq&\mbfP_{\mu}\left(\begin{aligned}\forall_{i< T,~\lfloor s\rfloor=i} ~x+\sigma_Q B_{s}+M^{(n,k)}_i\in[an^\alpha+2n^\rho,bn^\alpha-2n^\rho],
\\ x+\sigma_Q B_{T}+M^{(n,k)}_T\in[a'' n^\alpha+2n^\rho,b'' n^\alpha-2n^\rho]\end{aligned}\right)-c_3Tn^{-\rho\lambda_2}.
\end{eqnarray}

\emph{step 2: Using a Brownian motion to approximate $\{M_{t_{n,k}+i}-M_{t_{n,k}}\}_{i\leq T}$.}

We recall $\bar{t}_n:=t_n+n.$ Note that ${\mbfE}(M^2_1)=\sigma^2_A$. By Theorem \uppercase\expandafter{\romannumeral2}  and assumption (H2) we can find a standard Brownian motions $\hat{W}^{(n,k)}$ under $\mathbf{P}$ with $\hat{W}_0^{(n,k)}=0$~
such that
\begin{eqnarray}\label{3.3.3}~\forall n\in \bfN,~ \mathbf{P}\left(\max\limits_{i\leq T}|M_i^{(n,k)}- \hat{W}_{i\sigma_{_A}^2}^{(n,k)}|\geq n^{\rho}\right)\leq \frac{2(C^*\lambda_1)^{\lambda_1}T\mbfE(M_1^{\lambda_1})}{n^{\rho \lambda_1}}.\end{eqnarray}




Define
\begin{eqnarray}\label{3.3.4++++}E_n:=\cap_{k=0}^K\left\{\max\limits_{i\leq T}\left|M_i^{(n,k)}- \hat{W}_{{i\sigma_{_A}^2}}^{(n,k)}\right|\leq n^\rho\right\}\no,
\end{eqnarray}
\begin{eqnarray}\label{3.3.4}
\hat{E}_n:=\cap_{k=0}^K\left\{\sup_{s\in[i,i+1)} |\hat{W}_{{s\sigma_{_A}^2}}^{(n,k)}-\hat{W}_{{i\sigma_{_A}^2}}^{(n,k)}|\leq n^{\rho},~ \forall i\leq T\right\}.\end{eqnarray}
For any integer~$k\leq K-1,$ on the event $E_n\cap\hat{E}_n,$ we have
\begin{eqnarray}\label{3.3.7}
&&\mbfP_{\mu}\left(\begin{aligned}\forall_{i<T,~\lfloor s\rfloor=i} ~x+\sigma_Q B_{s}+M^{(n,k)}_i\in[an^\alpha+2n^\rho,bn^\alpha-2n^\rho],
\\ x+\sigma_Q B_{T}+M^{(n,k)}_T\in[a'' n^\alpha+2n^\rho,b'' n^\alpha-2n^\rho]\end{aligned}\right)\no
\\&=&\mbfP_{\mu}\left(\begin{aligned}\forall_{i<T,~\lfloor s\rfloor=i} ~x+\sigma_Q B_{s}+M^{(n,k)}_i-\hat{W}^{(n,k)}_{i\sigma^2_A}+\hat{W}^{(n,k)}_{i\sigma^2_A}\in[an^\alpha+2n^\rho,bn^\alpha-2n^\rho],
\\ x+\sigma_Q B_{T}+M^{(n,k)}_T-\hat{W}^{(n,k)}_{T\sigma^2_A}+\hat{W}^{(n,k)}_{T\sigma^2_A}\in[a'' n^\alpha+2n^\rho,b'' n^\alpha-2n^\rho]\end{aligned}\right)\no
\\&\geq&\mbfP_{\mu}\left(\begin{aligned}\forall_{i<T,~\lfloor s\rfloor=i} ~x+\sigma_Q B_{s}+\hat{W}^{(n,k)}_{i\sigma^2_A}\in[an^\alpha+3n^\rho,bn^\alpha-3n^\rho],
\\ x+\sigma_Q B_{T}+\hat{W}^{(n,k)}_{T\sigma^2_A}\in[a'' n^\alpha+3n^\rho,b'' n^\alpha-3n^\rho]\end{aligned}\right)\no
\\&\geq&\mbfP_{\mu}\left(\begin{aligned}\forall_{s<T} ~x+\sigma_Q B_{s}+\hat{W}^{(n,k)}_{s\sigma^2_A}\in[an^\alpha+4n^\rho,bn^\alpha-4n^\rho],
\\ x+\sigma_Q B_{T}+\hat{W}^{(n,k)}_{T\sigma^2_A}\in[a'' n^\alpha+3n^\rho,b'' n^\alpha-3n^\rho]\end{aligned}\right).
\end{eqnarray}
 According to \eqref{3.3.7}, \eqref{3.3.14} and the scaling property of Brownian motion, for any~$\varepsilon\in(0,\min\{b'-b'',b''-b''',a'''-a'',a''-a'\})$ and $n$ large enough satisfying $\varepsilon>4 n^{\rho-\alpha}$, on the event $F_n\cap\hat{F}_n\cap E_n\cap\hat{E}_n$ we have
\begin{eqnarray}\label{3.3.8+}
&&\inf_{x\in[a_0n^{\alpha},b_0n^{\alpha}]}\mbfP_\mu\left(\forall_{i\leq T}~ S_{t_{n,k}+i}\in[a n^\alpha, bn^\alpha], ~S_{t_{n,k+1}}\in[a''n^\alpha,b''n^\alpha] \Big|S_{t_{n,k}}=x\right)\no
\\&\geq&e^{X_{n,k}}-c_3Tn^{-\rho\lambda_2},
\end{eqnarray}
where $X_{n,k}:=\log\inf\limits_{x\in[a_0,b_0]}\mbfP_{\mu}\left(\begin{aligned}\forall_{s<T n^{-2\alpha}} ~x+\sigma_Q B_{s}+n^{-\alpha}\hat{W}^{(n,k)}_{s\sigma^2_An^{2\alpha}}\in[a+\varepsilon,b-\varepsilon],
\\ x+\sigma_Q B_{Tn^{-2\alpha}}+n^{-\alpha}\hat{W}^{(n,k)}_{T\sigma^2_An^{2\alpha}}\in[a''+\varepsilon,b''-\varepsilon]\end{aligned}\right), k<K.$

Note that the realization of $\hat{W}$ totally decide the quantity of $X_{n,k}$. That is to say, the influence on $X_{n,k}$ from $\mu$ is fully determined by $\hat{W}.$ Hence in the following proof we can write $X_{n,k}$ as $$\inf_{x\in[a_0,b_0]}\mbfP\left(\begin{aligned}\forall_{s<T n^{-2\alpha}} ~x+\sigma_Q B_{s}+n^{-\alpha}\hat{W}^{(n,k)}_{s\sigma^2_{_A}n^{2\alpha}}\in[a+\varepsilon,b-\varepsilon],
\\ x+\sigma_Q B_{Tn^{-2\alpha}}+n^{-\alpha}\hat{W}^{(n,k)}_{T\sigma^2_{_A}n^{2\alpha}}\in[a''+\varepsilon,b''-\varepsilon]\end{aligned}\bigg|\hat{W}\right).$$

For the last term in \eqref{3.3.2}, by similar discussion in \emph{step 1-2}, on the event $F_n\cap\hat{F}_n\cap E_n\cap\hat{E}_n$ we get
\begin{eqnarray}\label{3.3.16}
&&\inf\limits_{x\in [a'' n^\alpha,b'' n^\alpha]} \mbfP_\mu\left(\forall_{i\leq \bar{t}_n-t_{n,K}}~ S_{t_{n,K}+i}\in[an^\alpha, b n^\alpha], S_{\bar{t}_n}\in[a'n^\alpha,b'n^\alpha] \big|S_{t_{n,K}}=x\right)\no
\\&\geq&\inf\limits_{x\in [a'' n^\alpha,b'' n^\alpha]} \mbfP_\mu\left(\forall_{i\leq T}~ S_{t_{n,K}+i}\in[a'n^\alpha, b'n^\alpha]\big|S_{t_{n,K}}=x\right)\no
\\&\geq&\inf_{x\in[a'',b'']}\mbfP\left(\forall_{s<T n^{-2\alpha}} ~x+\sigma_Q B_{s}+n^{-\alpha}\hat{W}^{(n,K)}_{s\sigma_{_A}^2 n^{2\alpha}}\in[a'+\varepsilon,b'-\varepsilon]\big|\hat{W}\right)-c_3Tn^{-\rho\lambda_2}\no
\\&\geq&e^{X_{n,K}}-c_3Tn^{-\rho\lambda_2},\end{eqnarray}
where $X_{n,K}:=\log\inf_{x\in[a'',b'']}\mbfP\left(\forall_{s<D} ~x+\sigma_Q B_{s}+n^{-\alpha}\hat{W}^{(n,K)}_{s\sigma_{_A}^2n^{2\alpha}}\in[a'+\varepsilon,b'-\varepsilon]\big|\hat{W}\right).$

\emph{step 3. Define
\begin{eqnarray}\label{3.3.18}G_{n}:=\left\{\min_{k\leq K} \frac{e^{X_{n,k}}}{2}\geq c_3Tn^{-\rho\lambda_2}\right\}~~{\rm and}~~ H_n:=E_n\cap \hat{E}_n\cap F_n\cap\hat{F_n}\cap G_n. \end{eqnarray}
We find the lower bound on the event $H_n.$}

By Markov property we have
\begin{eqnarray}\label{3.3.19}X_{n,k}\geq X_{n,k,D-1}+X_{n,k,1},~~~\forall k\leq K-1.\end{eqnarray} where 
 $$X_{n,k,D-1}:=\log\inf\limits_{x\in [a_0,b_0]}\mbfP\left(\begin{aligned}\forall_{s\leq D-1}~ x+\sigma_{_Q} B_{s}+n^{-\alpha}\hat{W}^{(n,k)}_{s\sigma_{_A}^2n^{2\alpha}}\in[a+\varepsilon, b-\varepsilon],
\\ x+\sigma_{_Q} B_{D-1}+ n^{-\alpha}\hat{W}^{(n,k)}_{_{(D-1)\sigma_{_A}^2n^{2\alpha}}}\in[a''',b''']\end{aligned}\Bigg|\hat{W}\right),$$
$$ X_{n,k,1}:=\log\inf\limits_{x\in [a''',b''']}\mbfP\left(\forall_{s\leq 1}~x+
\sigma_{_Q}B_{s}+\frac{\hat{W}^{(n,k)}_{_{(D-1+s)\sigma_{_A}^2n^{2\alpha}}}-
\hat{W}^{(n,k)}_{_{(D-1)\sigma_{_A}^2n^{2\alpha}}}}{n^{\alpha}}\in[a''+\varepsilon, b''-\varepsilon]\big|\hat{W}\right).$$
Combining~\eqref{3.3.2} with~\eqref{3.3.16}-\eqref{3.3.19}, we get
\begin{eqnarray}\label{3.3.20}
&&\frac{\log\inf\limits_{x\in[a_0 n^{\alpha}, b_0 n^{\alpha}]}\mbfP_\mu
\Big(\forall_{t_n\leq i\leq \bar{t}_n}~ S_{i}\in[an^{\alpha},bn^{\alpha}],~S_{\bar{t}_n}\in[a'n^{\alpha}, b'n^{\alpha}] \Big|S_{t_n}=x\Big)}{n^{1-2\alpha}}\nonumber
\\&\geq &\frac{\log\left[\1_{H_n}\inf\limits_{x\in[a_0 n^{\alpha}, b_0 n^{\alpha}]}\mbfP_\mu
\Big(\forall_{t_n\leq i\leq \bar{t}_n}~ S_{i}\in[an^{\alpha},bn^{\alpha}],~S_{\bar{t}_n}\in[a'n^{\alpha}, b'n^{\alpha}] \big|S_{t_n}=x\Big)\right]}{n^{1-2\alpha}}\nonumber
\\&\geq &\frac{\log\1_{H_n}+\sum_{k=0}^{K-1}(-\log 2+X_{n,k,D-1}+X_{n,k,1})-\log 2+X_{n,K}}{K}\frac{K}{n^{1-2\alpha}}.
\end{eqnarray}
Let $W$ be a standard Brownian motion which is independent of $B.$ By the scaling property of Brownian motion we know that for 
any~$n,~k,$~random variables $X_{n,k,D-1},$~$X_{n,k,1},~X_{n,K}$ have the same distributions as~$~Y_{_{D-1}},~Y_{1},~Y_{_D}$ respectively,~where
$$Y_{_{D-1}}:=\log\inf\limits_{x\in [a_0,b_0]}\mbfP\left(\begin{aligned}\forall_{s\leq D-1}~ x+\sigma_Q B_{s}+ \sigma_A W_s\in[a+\varepsilon, b-\varepsilon],\\
 x+\sigma_Q B_{D-1}+\sigma_A W_{D-1}\in[a''',b''']\end{aligned}\Bigg|W\right),$$
$$Y_{1}:=\log\inf\limits_{x\in[a''',b''']} \mbfP\big(\forall_{0\leq s\leq 1}~ x+\sigma_Q B_{s}+\sigma_A W_s\in [a''+\varepsilon,b''-\varepsilon] |W),$$
$$Y_{_D}:= \log\inf\limits_{x\in [a'',b'']} \mbfP\Big(\forall_{s\leq D}~ x+\sigma_Q B_{s}+\sigma_A W_s\in[a'+\varepsilon, b'-\varepsilon]|W\Big).$$  Theorem \uppercase\expandafter{\romannumeral1} \eqref{bmt++} shows that
\begin{eqnarray}\label{3.3.20}\mbfE(|Y_{_{D-1}}|^j)+\mbfE(|Y_{1}|^j)+\mbfE(|Y_{D}|^j)<+\infty,~~ \forall j\in\bfN.\end{eqnarray} Recalling that $D\in\bfN^+, ~T:=\lfloor Dn^{2\alpha}\rfloor,~ K:=\left\lfloor \frac{n}{T} \right\rfloor.$ According to \eqref{3.3.20} and~Borel-Cantelli 0-1 law we can obtain
\begin{eqnarray}\label{3.3.21}
\lim\limits_{n\rightarrow+\infty}\frac{\sum_{k=0}^{K-1}X_{n,k,D-1}}{K}=\mathbf{E}(Y_{D-1}),~{\rm \mathbf{P}-a.s..}\end{eqnarray}
\begin{eqnarray}\label{3.3.22}~\lim\limits_{n\rightarrow+\infty}\frac{\sum_{k=0}^{K-1}X_{n,k,1}}{K}=\mathbf{E}(Y_{1}),~
~\lim\limits_{n\rightarrow+\infty}\frac{X_{n,K}}{K}=0,~{\rm \mathbf{P}-a.s.}.
\end{eqnarray}
Since the proof of the three limits in \eqref{3.3.21} and \eqref{3.3.22} are similar, we only show the proof of \eqref{3.3.21}.
~Denote~$\hat{X}_{n,k,D-1}:=X_{n,k,D-1}-\mathbf{E}(Y_{D-1}).$ Choose an even~$m$ which satisfies~$\frac{(1-2\alpha)m}{2}>1,$ then for any~$\epsilon>0,$ we have
\begin{eqnarray}\label{3.3.23}
\mathbf{P}\left(\left|\frac{\sum_{k=0}^{K-1}X_{n,k,D-1}}{K}-\mathbf{E}(Y_{D-1})\right|>\epsilon\right)&=&\mathbf{P}\left(\left|\frac{\sum_{k=0}^{K-1}\hat{X}_{n,k,D-1}}{K}\right|>\epsilon\right)\nonumber
\\&\leq&\frac{\mathbf{E}\left[(\sum_{k=0}^{K-1}\hat{X}_{n,k,D-1})^m\right]}{K^m\epsilon^{m}}.
\end{eqnarray}
Now we will show \begin{eqnarray}\label{3.3.23+}\mathbf{E}\left[\left(\sum_{k=0}^{K-1}\hat{X}_{n,k,D-1}\right)^m\right]\leq O(K^{\frac{m}{2}}).\end{eqnarray}

To do this, the most important basis is the fact that for any fixed~$n$, $\{\hat{X}_{n,k,D-1}\}_{k\leq K-1}$ is an i.i.d. centered random sequence with ~
finite $j-$th moment, $\forall j\in\bfN$. Hence when we expand~$(\sum_{k=0}^{K-1}\hat{X}_{n,k,D-1})^m$,~the expectations of the terms which contain the first moment are $0.$ For the rest terms, if we call the form like ``$\Pi^{l}_{i=1}\hat{X}^{\iota_i}_{n,k_i,D-1}, k_1<k_2<\cdots<k_l$'' the $l-$ elements term, then when $l>\frac{m}{2},$ there is no $l-$ elements term in the rest terms because every term in the rest terms satisfies $\min_{1\leq i\leq l}\iota_i\geq 2$ and $\sum_{i=1}^{l}\iota_i=m.$ Moreover, the number of the $l-$ elements term is at most $\mathcal{C}^{l}_{K} l^{m},$ where $\mathcal{C}^{\cdot}_{\cdot}$ represents the combination number, i.e.,~$\mathcal{C}^{l}_{K}=\frac{K!}{l!(K-l)!}.$
On the other hand, we observe that the expectation of any term in the expansion of $(\sum_{i=0}^{K-1}\hat{X}_{n,k,D-1})^m$ will not larger than the constant $C_{Y}:=\max\left([\max_{j\leq m}\mbfE(|Y_{_{D-1}}-\mbfE Y_{_{D-1}}|^j)]^m,1\right).$ Moreover, \eqref{3.3.20} means that $C_Y<+\infty.$ In conclusion, for $K$ large enough such that $K>m,$ we have
\begin{eqnarray}\label{3.3.23++}\mathbf{E}\left[(\sum_{i=0}^{K-1}\hat{X}_{n,k,D-1})^m\right]&\leq& \sum_{l=1}^{\frac{m}{2}}\mathcal{C}^{l}_{K} l^{m}C_Y\leq \frac{m}{2}\mathcal{C}^{\frac{m}{2}}_{K}\left(\frac{m}{2}\right)^{m}C_Y\leq K^{\frac{m}{2}}\left(\frac{m}{2}\right)^mC_Y.\no
\end{eqnarray}
This finishes the proof of \eqref{3.3.23+}\footnote{We note that the quantities $l^{m}$ and $C_Y$ as upper bounds can be improved but it is not relevant for our proof.}. We recall the relationships $K=O(n^{1-2\alpha})$ and $\frac{(1-2\alpha)m}{2}>1,$ which combined with \eqref{3.3.23} and \eqref{3.3.23+} imply that
$$\sum_{n=1}^{+\infty}\mathbf{P}\left(\left|\frac{\sum_{k=0}^{K-1}X_{n,k,D-1}}{K}-\mathbf{E}(Y_{D-1})\right|>\epsilon\right)\leq
 \sum_{n=1}^{+\infty}\frac{O(K^{m/2})}{K^m\epsilon^{m}}<+\infty.$$ 
Therefore, we get \eqref{3.3.21} by~Borel-Cantelli 0-1 law.

Combining the above discussions in \emph{step 3} ~we finally get
\begin{eqnarray}\label{3.3.24}
&&\varliminf\limits_{n\rightarrow+\infty}\frac{\log\inf\limits_{x\in[a_0 n^{\alpha}, b_0 n^{\alpha}]}\mbfP_\mu
\Big(\forall_{t_n\leq i\leq \bar{t}_n}~ S_{i}\in[an^{\alpha},bn^{\alpha}],~S_{\bar{t}_n}\in[a'n^{\alpha}, b'n^{\alpha}] \Big|S_{t_n}=x\Big)}{n^{1-2\alpha}}\nonumber
\\&\geq&\frac{-\log 2+\mathbf{E}(Y_{D-1})+\mathbf{E}(Y_{1})}{D}+\varliminf\limits_{n\rightarrow+\infty}\frac{\log\1_{H_n}}{D}.
\end{eqnarray}

\emph{step 4. Showing that $\lim\limits_{n\rightarrow+\infty}\1_{H_n}=1, ~{\rm \mathbf{P}-a.s.}$ if $\lambda_0=2\alpha^{-1};$ while $\lim\limits_{n\rightarrow+\infty}\1_{H_n}=1,$ in Probability~${\mathbf{P}}$ if $\lambda_0=\alpha^{-1}.$}

We recall $H_n:=E_n\cap \hat{E}_n\cap F_n\cap \hat{F}_n\cap G_n.$

Thanks to~\eqref{3.3.3}, it is true that $$\mbfP(E^c_n)\leq (K+1)2(C^*\lambda_1)^{\lambda_1}\mbfE(M_1^{\lambda_1})Tn^{-\rho \lambda_1}\leq4(C^*\lambda_1)^{\lambda_1}\mbfE(M^{\lambda_1}_1)n^{1-\rho\lambda_1}.$$ Hence we have $\lim\limits_{n\rightarrow+\infty}\1_{E_n}=1, ~{\rm in~Probability~ \mathbf{P}}$ if $\lambda_0=\alpha^{-1}$. Moreover, by Borel-Cantelli 0-1 law we can see $\lim\limits_{n\rightarrow+\infty}\1_{E_n}=1, ~{\rm \mathbf{P}-a.s.}$ if $\lambda_0=2\alpha^{-1}.$

By Cs\"{o}rg\H{o} and R\'{e}v\'{e}sz~estimation (\cite[Lemma~1]{CR1979}), we know that there exist two positive constants~$c_4,c_5$ such that
$$\mathbf{P}(\hat{E}^c_n)\leq (K+1)T\mathbf{P}\left(\sup_{s\in[0,1]}|W_s|>n^{\rho}\right)\leq c_4 n e^{-c_5n^{2\rho}}, ~\forall n\in\bfN.$$
By Borel-Cantelli 0-1 law and the fact~$\sum\limits_{n=1}^{+\infty}c_1 n e^{-c_4n^{c_5\rho}}<+\infty,$ we get~$\lim\limits_{n\rightarrow+\infty}\1_{\hat{E}_n}=1, ~{\rm \mathbf{P}-a.s..}$

Next we estimate $\hat{F}_n.$ Let $\bar{\zeta}_{i}:={\zeta}_{i}-i\mbfE({\zeta}_{1}).$  By Theorem \uppercase\expandafter{\romannumeral2} we can construct a standard Brownian motion $W^{(\zeta)}$ such that $\mbfP\left(\max_{i\leq T}|\bar{\zeta}_{i}-W^{(\zeta)}_{i\mbfE(\bar{\zeta}^2_1)}|>C^*\lambda_3n^{2\rho}\right)\leq 2T\mbfE(\bar{\zeta}_{1}^{\lambda_3})n^{-2\rho\lambda_3}.$ Note that $\rho<\alpha$ and the Cs\"{o}rg\H{o} and R\'{e}v\'{e}sz~estimation mentioned above. Then for large enough $n$ we have
\begin{eqnarray}\mathbf{P}\Big(\zeta_{n,0}\geq 2\mbfE(\zeta_1)T\Big)&\leq& \mathbf{P}\Big(|\bar{\zeta}_T|\geq \mbfE(\zeta_1)T, W^{(\zeta)}_{T\mbfE(\bar{\zeta}^2_1)}\leq n^{2\rho}\Big)+\mathbf{P}\Big(W^{(\zeta)}_{T\mbfE(\bar{\zeta}^2_1)}> n^{2\rho}\Big)\no
\\&\leq& \mbfP\left(\max_{i\leq T}|\bar{\zeta}_{i}-W^{(\zeta)}_{i\mbfE(\bar{\zeta}^2_1)}|>C^*\lambda_3n^{2\rho}\right) +\mathbf{P}\Big(W^{(\zeta)}_{T\mbfE(\bar{\zeta}^2_1)}> n^{2\rho}\Big)\no
\\&\leq& 3\mbfE(\bar{\zeta}_{1}^{\lambda_3})Tn^{-2\rho\lambda_3}.\no\end{eqnarray}
By the same method we get $\mathbf{P}\Big(\zeta_{n,k}\geq 2\mbfE(\zeta_1)T\Big)\leq 3\mbfE(\bar{\zeta}_{1}^{\lambda_3})Tn^{-2\rho\lambda_3}, \forall k\leq K.$ Therefore, $\mbfP(\hat{F}^c_n)\leq (K+1)3\mbfE(\bar{\zeta}_{1}^{\lambda_3})Tn^{-2\rho\lambda_3}\leq 4\mbfE(\bar{\zeta}_{1}^{\lambda_3})n^{1-2\rho\lambda_3},$ which means that $\lim\limits_{n\rightarrow+\infty}\1_{\hat{F}_n}=1, ~{\rm \mathbf{P}-a.s.}$ when $\lambda_0=2\alpha^{-1},$ while $\lim\limits_{n\rightarrow+\infty}\1_{\hat{F}_n}=1,$ in Probability~${\mathbf{P}}$ when $\lambda_0=\alpha^{-1}.$

We turn to $F_n.$ The approach is similar to what we do for the $\hat{F}_n.$ Denote $\bar{\Gamma}_{i}:={\Gamma}_{i}-i\sigma^2_Q$. By Assumption (H3) and the Jensen's inequality we have $\mbfE(\Gamma_1^{\lambda_2\lambda_3/2})<+\infty$ thus $\mbfE(\Gamma_1^{\lambda_3})<+\infty.$ By Theorem \uppercase\expandafter{\romannumeral2} we can construct a standard Brownian motion $W^{(\Gamma)}$ such that
$\mbfP(\max_{i\leq T}|\bar{\Gamma}_{i}-W^{(\Gamma)}_{i\sigma_{Q}^2}|\geq C^*\lambda_3n^{\alpha+\varsigma})\leq 2\mbfE(\bar{\Gamma}_1^{\lambda_3})T n^{-(\alpha+\varsigma) \lambda_3}.$ Hence for large enough $n$ we have
\begin{eqnarray}&&\mathbf{P}\Big(\max_{0\leq i \leq T}|\Gamma_i-i\sigma^2_Q|\geq 2C^*\lambda_3n^{\alpha+\varsigma}\Big)\no
\\&\leq& \mathbf{P}\Big(\max_{0\leq i \leq T}|\bar{\Gamma}_i|\geq 2C^*\lambda_3n^{\alpha+\varsigma}, \max_{0\leq i \leq T}|W^{(\Gamma)}_{i\sigma_{Q}^2}|\leq C^*\lambda_3n^{\alpha+\varsigma}\Big)+\mathbf{P}\Big(\max_{0\leq i \leq T}|W^{(\Gamma)}_{i\sigma_{Q}^2}|> C^*\lambda_3n^{\alpha+\varsigma}\Big)\no
\\&\leq& \mbfP(\max_{0\leq i\leq T}|\bar{\Gamma}_{i}-W^{(\Gamma)}_{i\sigma_{Q}^2}|\geq C^*\lambda_3n^{\alpha+\varsigma})+\mathbf{P}\Big(\max_{0\leq i \leq T}|W^{(\Gamma)}_{i\sigma^2_{Q}}|> C^*\lambda_3n^{\alpha+\varsigma}\Big)\no
\\&\leq& 3\mbfE(\bar{\Gamma}^{\lambda_3}_{1})Tn^{-(\alpha+\varsigma)\lambda_3},\no\end{eqnarray}
which means that $\lim\limits_{n\rightarrow+\infty}\1_{\hat{F}_n}=1, ~{\rm \mathbf{P}-a.s.}$ if $\lambda_0=2\alpha^{-1},$ while $\lim\limits_{n\rightarrow+\infty}\1_{\hat{F}_n}=1,$ in Probability~${\mathbf{P}}$ if $\lambda_0=\alpha^{-1}.$

At last we consider~$G_n.$ We should note that $\mathbf{E}(e^{c_6|Y_{D-1}|}+e^{c_6|Y_{1}|}+e^{c_6|Y_{D}|})<+\infty$ for any $c_6>0$ because of \eqref{bmt+++}. Recall the definition of $G_n$ and the relationship $2\alpha<\rho\lambda_2$, then for any~$\epsilon>0$ and large enough $n$, we can find a constant $\nu>0$ such that
 \begin{eqnarray}\mathbf{P}(|\1_{G_{n}}-1|\geq \epsilon)&=&\mathbf{P}(G^c_{n})\no
 \\&\leq &K\mathbf{P}(-\log 2-X_{n,0,D-1}-X_{n,0,1}>3\nu\log n)+\mathbf{P}(-\log 2-X_{n,K}>2\nu\log n)\nonumber
\\&\leq &K[\mathbf{P}(|Y_{D-1}|> \nu\log n)+\mathbf{P}(|Y_{1}|> \nu\log n)] + \mathbf{P}(|Y_{D}|>\nu\log n).\nonumber
\\&\leq &2Kn^{-c_7\nu}[\mathbf{E}(e^{c_7|Y_{D-1}|}+e^{c_7|Y_{1}|}+e^{c_7|Y_{D}|})],
\end{eqnarray}
where we choose $c_7>\nu^{-1}(2-2\alpha).$
Hence we have $\sum_{n=1}^{+\infty}\mbfP(|\1_{G_{n}}-1|\geq \epsilon)<+\infty.$
By~Borel-Cantelli 0-1 law we see~$\lim\limits_{n\rightarrow+\infty}\1_{G_n}=1, ~{\rm \mathbf{P}-a.s.}.$


In conclusion, finally we get
\begin{eqnarray}\label{3.3.25}\lim\limits_{n\rightarrow+\infty}\1_{H_n}=1, ~{\rm \mathbf{P}-a.s.} ~\text{if}~ \lambda_0=2\alpha^{-1}; \text{~in ~Probability}~{\rm \mbfP} ~\text{if}~ \lambda_0=\alpha^{-1}.\end{eqnarray}

 Let us go back to~\eqref{3.3.24}. As~$D\ra +\infty,$ the $L^1$ convergence in Theorem \uppercase\expandafter{\romannumeral1} implies that $$\lim\limits_{D\rightarrow+\infty}\frac{\mathbf{E}(Y_{D-1})}{D}=\frac{-\tilde{\gamma}}{(b-a-2\varepsilon)^2}.$$
Combining with~\eqref{3.3.25} and the fact~$\mathbf{E}(Y_{1})<+\infty$ , we get the lower bound~\eqref{2-}~by taking $\varepsilon\ra 0$. This completes the proof of the lower bound. \qed
\subsection{The upper bound of Theorem 2.2}

~~The approach to the upper bound is analogous to the proof of the lower bound. Hence we only present a sketch of this proof.

{\bf Proof of Theorem \ref{mog}: the upper bound \eqref{1+}}

Here we let $D\in\bfN^+, T:=\lceil Dn^{2\alpha}\rceil, K:=\lfloor\frac{n}{T}\rfloor.$
~Adjusting the definition of $G_n$ to $$G_n:=\Big\{\min_{k\leq K-1}\sup_{x\in\bfR}\mbfP\Big(\forall_{s\leq D} x+\sigma_Q B_{s}+n^{-\alpha}\hat{W}^{(n,k)}_{s\sigma^2_An^{2\alpha}}\in[a-\varepsilon, b+\varepsilon]\big|\hat{W}\Big)\geq c_3Tn^{-\rho\lambda_2}\Big\}.$$
Other notations and definitions in Section 3.2 do not need to be changed. By Markov property we see
\beqnn
&&\log\sup\limits_{x\in \bfR}\mbfP_\mu
(\forall_{t_n\leq i\leq \bar{t}_n} an^\alpha\leq S_{i}\leq bn^\alpha|S_{t_n}=x)
\\&\leq&\mathbf{1}_{H_n}\log\sup\limits_{x\in \bfR}\mbfP_{\mu}
(\forall_{t_n\leq i\leq \bar{t}_n} ~x+U^{(n,k)}_{i}+M^{(n,k)}_{i} \in[an^\alpha,bn^\alpha])
\\&\leq&\mathbf{1}_{H_n}\log\left\{\prod\limits_{k=0}\limits^{K-1}\sup\limits_{x\in\bfR} \mbfP_\mu(\forall_{i\leq T} ~x+U^{(n,k)}_{i}+M^{(n,k)}_{i}\in[an^\alpha, bn^\alpha])\right\}.
\eeqnn
Taking the similar strategy as \emph{step 1} in the proof of the lower bound, for any $k=0,1,...K-1,$ on the event $H_n$ we have
\beqnn
&&\mbfP_\mu(\forall_{i\leq T} ~x+U^{(n,k)}_{i}+M^{(n,k)}_{i}\in[an^\alpha, bn^\alpha])
\\&\leq& \mbfP_{\mu}(\forall_{i\leq T}~x+B_{\Gamma^{(n,k)}_{i}}+M^{(n,k)}_{i}\in[an^\alpha-n^\rho, bn^\alpha+n^\rho])+\mbfP_\mu(\mathcal{V}^{c}_{n,k})
\\&\leq& \mbfP_{\mu}(\forall_{i\leq T,~\lfloor s\rfloor=i}~x+B_{s\sigma^2_Q}+M^{(n,k)}_{i}\in[an^\alpha-2n^\rho, bn^\alpha+2n^\rho])
+\mbfP_\mu(\mathcal{V}^{c}_{n,k})+\mbfP_\mu(\Theta^c_{n})
\\&\leq&\mbfP_{\mu}(\forall_{i\leq T,~\lfloor s\rfloor=i}~x+B_{s\sigma^2_Q}+M^{(n,k)}_{i}\in[an^\alpha-2n^\rho, bn^\alpha+2n^\rho])+c_3Tn^{-\rho\lambda_2}.
\eeqnn
Taking the similar strategy as \emph{step 2} in the proof of the lower bound, for any $k=0,1,...K-1,$ on the event $H_n$ we have
\beqnn
&&\mbfP_{\mu}\left(\forall_{i\leq T,~\lfloor s\rfloor=i}~x+B_{s\sigma^2_Q}+M^{(n,k)}_{i}\in[an^\alpha-2n^\rho, bn^\alpha+2n^\rho]\right)
\\&\leq& \mbfP\left(\forall_{i\leq T,~\lfloor s\rfloor=i}~x+B_{s\sigma^2_Q}+\hat{W}^{(n,k)}_{i\sigma^2_A}\in[an^\alpha-3n^\rho, bn^\alpha+3n^\rho]\big|\hat{W}^{(n,k)}\right)
\\&\leq& \mbfP\left(\forall_{s\leq T}~~x+B_{s\sigma^2_Q}+\hat{W}^{(n,k)}_{s\sigma^2_A}\in[an^\alpha-4n^\rho, bn^\alpha+4n^\rho]\big|\hat{W}^{(n,k)}\right).
\eeqnn
Note that $\frac{T}{n^{2\alpha}}\geq D,$ by the scaling property of Brownian motion and the definition of $G_n,$ we can see
\beqnn&&\log\sup\limits_{x\in\bfR}\mbfP_\mu
(\forall_{t_n\leq i\leq \bar{t}_n}~ an^\alpha\leq S_{i}\leq bn^\alpha|S_{t_n}=x)
\\&\leq&\mathbf{1}_{H_n}\left\{\sum\limits_{k=0}\limits^{K-1}\log\Big[\sup\limits_{x\in\bfR} \mbfP\Big(\forall_{s\leq D}~x+\sigma_{_Q} B_{s}+n^{-\alpha}\hat{W}^{(n,k)}_{s\sigma^2_An^{2\alpha}}\in[a-\varepsilon, b+\varepsilon]|\hat{W}^{(n,k)}\Big)+c_3Tn^{-\rho\lambda_2}\Big]\right\}
\\&\leq&\mathbf{1}_{H_n}\left\{\sum\limits_{k=0}\limits^{K-1}\log\Big[2\sup\limits_{x\in\bfR} \mbfP\Big(\forall_{s\leq D}~x+\sigma_{_Q} B_{s}+n^{-\alpha}\hat{W}^{(n,k)}_{s\sigma^2_An^{2\alpha}}\in[a-\varepsilon, b+\varepsilon]|\hat{W}^{(n,k)}\Big)\Big]\right\}.
\eeqnn
After the similar discussions as the \emph{step 3-4} in the proof of the lower bound, finally we obtain
\beqnn&&\varlimsup_{n\rightarrow+\infty}\frac{\log\sup\limits_{x\in\bfR}\mbfP_\mu
\left(\forall_{t_n\leq i\leq \bar{t}_n} an^\alpha\leq S_{i}\leq bn^\alpha|S_{t_n}=x\right)}{n^{1-2\alpha}}
\\&\leq&\frac{\log2+\mbfE\big[\log\sup\limits_{x\in\bfR}\mbfP\big(\forall_{s\leq D}~ \sigma_Q B_{s}+\sigma_{A}W_s\in[a-\varepsilon, b+\varepsilon]|W\big)\big]}{D}.
\eeqnn
At last, let $D\rightarrow+\infty,$ $\varepsilon\rightarrow 0,$ then we can get the upper bound \eqref{upper} by applying the $L^1$ convergence in Theorem \uppercase\expandafter{\romannumeral1}.\qed

\subsection{ Proof of Corollary \ref{mogc1}}

Now we prove the Corollary 2.1. The proof mainly depends on the Markov property and the continuity of the functions $g$ and $h$.
Here we only prove the lower bound. The upper bound is similar and easier than the lower bound so we left it to readers.
Choose $m\in\bfN$ and let $d_n:=\lfloor\frac{n}{m}\rfloor, \tilde{t}_{n,k}:=t_n+kd_n.$ Without loss of generality, in this proof
we assume that $a_0<a'<0<b'<b_0.$ For any $a,b$ satisfying $a<a_0, b>b_0,$ we choose $v>0$ such that $(\frac{m+v}{m})^{\alpha}<\min\{\frac{a}{a_0},\frac{b}{b_0}\}.$ Let $n$ large enough such that $n\leq (m+v)d_n,$ then
\beqnn
Q_{d_n,k}&:=&\inf\limits_{x\in[a_0 n^{\alpha}, b_0 n^{\alpha}]}\mbfP_\mu
\left(\forall_{0\leq i\leq d_n}~ S_{\tilde{t}_{n,k}+i}\in[an^\alpha, bn^\alpha], S_{\tilde{t}_{n,k+1}}\in[a'n^\alpha, b'n^\alpha]\Big|S_{\tilde{t}_{n,k}}=x\right)
\\&\geq&\inf\limits_{x\in[a_0(m+v)^{\alpha} d_n^{\alpha}, b_0(m+v)^{\alpha} d_n^{\alpha}]}\mbfP_\mu
\left(\begin{aligned}\forall_{0\leq i\leq d_n}~ S_{\tilde{t}_{n,k}+i}\in[(am^\alpha)d_n^{\alpha}, (bm^\alpha)d_n^{\alpha}],\\
S_{\tilde{t}_{n,k+1}}\in[(a'm^\alpha)d_n^{\alpha}, (b'm^\alpha)d_n^{\alpha}]\end{aligned}\Bigg|S_{\tilde{t}_{n,k}}=x\right).
\eeqnn
Therefore, by Theorem 2.2 we have \beqlb\label{mogc1001}\varliminf_{n\rightarrow +\infty}\frac{\log Q_{d_n,k}}{n^{1-2\alpha}}
\geq-m^{2\alpha-1}\frac{\tilde{\gamma}}{(bm^{\alpha}-am^{\alpha})^2}
=-\frac{1}{m}\frac{\tilde{\gamma}}{(b-a)^2},~~~{\rm \mathbf{P}-a.s..}
\eeqlb
Let $\bar{\varepsilon}:=\min\left\{\frac{\inf_{s\in[0,1]}(h(s)-g(s))}{6}, \frac{b_0-a_0}{2},\frac{b'-a'}{2}, h(0)-b_0, a_0-g(0),h(1)-b', a'-g(1)\right\}.$ Choosing an $\varepsilon\in (0,\bar{\varepsilon}),$ then there exists a function $l$ defined on $[0,1]$ such that $$l(0)\in(a_0+\varepsilon,b_0-\varepsilon),~l(1)\in(a'+\varepsilon,b'-\varepsilon),~l(t)\in(g(t)+3\varepsilon,h(t)-3\varepsilon),~~ \forall t\in(0,1).$$  Choose $m\in\bfN$ large enough such that $$\sup\limits_{0\leq|t-s|\leq2m^{-1}}\{|g(t)-g(s)|+|h(t)-h(s)|\}\leq \varepsilon. $$
Denote $J_{k}:=\big[(l(k/m)-\varepsilon)n^{\alpha},(l(k/m)+\varepsilon)n^{\alpha}\big]$ and $I_{k,j}:=\left[g\Big(\frac{k d_n+j}{n}\Big)n^\alpha, h\Big(\frac{k d_n+j}{n}\Big)n^\alpha\right].$
By recalling $\bar{t}_n:=t_n+n$ and $\tilde{t}_{n,k}:=t_n+kd_n,$ we observe that
\beqlb\label{mogc1002}
&&\inf\limits_{x\in [a_0 n^\alpha, b_0 n^\alpha]}\mbfP_{\mu}\left(\forall_{0\leq j\leq n} S_{t_n+j}\in I_{0,j},~S_{\bar{t}_n}\in[a'n^\alpha,b'n^\alpha]\Bigg|S_{t_n}=x\right)\no
\\&\geq& \inf\limits_{x\in [a_0 n^\alpha, b_0 n^\alpha]}\mbfP_{\mu}\left(\forall_{0\leq j\leq d_{n}} S_{t_n+j}\in I_{0,j}, S_{\tilde{t}_{n,1}}\in J_1\Big|S_{t_n}=x\right)\no
\\&\times&\prod\limits_{k=1}\limits^{m-1}\inf\limits_{x\in J_k}\mbfP_{\mu}\left(\forall_{0\leq j\leq d_{n}}\begin{aligned} S_{\tilde{t}_{n,k}+j}\in I_{k,j} , S_{\tilde{t}_{n,k+1}}\in J_{k+1}\end{aligned}|S_{\tilde{t}_{n,k}}=x\right)\no
\\&\times&\inf\limits_{x\in J_{m}}\mbfP_{\mu}(\forall_{0\leq j\leq d_{n}} S_{\tilde{t}_{n,m}+j}\in[a'n^\alpha ,b'n^\alpha]|S_{\tilde{t}_{n,m}}=x\Big),
\eeqlb
Let $a\vee b:=\max\{a,b\}$ and $a\wedge b:=\min\{a,b\}.$ One checks that
$$\left(0\vee\frac{k-1}{m}\right)\leq\frac{k d_n}{n}\leq\frac{(k+1)d_n}{n}\leq\left(1\wedge\frac{k+1}{m}\right)$$
as long as $n>m^2.$ Denote $$\overline{g}_{k,m}:=\sup_{x\in\left[0\vee\frac{k-1}{m},1\wedge\frac{k+1}{m}\right]}g(x),~~~\underline{h}_{k,m}:=\inf_{x\in\left[0\vee\frac{k-1}{m},1\wedge\frac{k+1}{m}\right]}h(x).$$
According to the choices of~$\varepsilon$ and the function~$l$, for each $k\in\{0,1,...,m-1\},$ we can see that~$$\min\left\{l\left(\frac{k}{m}\right),l\left(\frac{k+1}{m}\right)\right\}-\varepsilon>\overline{g}_{k,m},~~\max\left\{l\left(\frac{k}{m}\right),l\left(\frac{k+1}{m}\right)\right\}+\varepsilon<\underline{h}_{k,m}.~~$$ 
Note that the RWre is space-homogeneous, hence we have
\beqnn
&&\inf\limits_{x\in J_k}\mbfP_{\mu}\left(\forall_{0\leq j\leq d_{n}}\begin{aligned} S_{\tilde{t}_{n,k}+j}\in I_{k,j} , S_{\tilde{t}_{n,k+1}}\in J_{k+1}\end{aligned}|S_{\tilde{t}_{n,k}}=x\right)
\\&=&\inf\limits_{x+l(\frac{k+1}{m})n^{\alpha}\in J_k}\mbfP_{\mu}\left(\forall_{0\leq j\leq d_{n}}\begin{aligned} S_{\tilde{t}_{n,k}+j}+l\left(\frac{k+1}{m}\right)n^{\alpha}\in I_{k,j} , S_{\tilde{t}_{n,k+1}}\in [-\varepsilon n^{\alpha},\varepsilon n^{\alpha}]\end{aligned}|S_{\tilde{t}_{n,k}}=x\right).\eeqnn
Combining with the above analysis and applying~\eqref{mogc1001} to each term in~\eqref{mogc1002},
finally we get
\beqnn
&&\varliminf_{n\rightarrow +\infty}\frac{\log\inf\limits_{x\in [a_0 n^\alpha, b_0 n^\alpha]}\mbfP_{\mu}\left(\forall_{0\leq j\leq n} S_{t_n+j}\in I_{0,j},~S_{\bar{t}_n}\in[a'n^\alpha,b'n^\alpha]|S_{t_n}=x\right)}{n^{1-2\alpha}}
\\&\geq& \frac{1}{m}\left[\sum_{k=0}^{m-1}\frac{\tilde{\gamma}}{(\underline{h}_{k,m}-\overline{g}_{k,m})^2}
+\frac{\tilde{\gamma}}{(b'-a')^2}\right].
\eeqnn
Hence we complete the proof by taking $m\rightarrow +\infty.$ \qed

\section{Appendix: Proof of \eqref{bmt+++}}
In this proof we keep all the notations which is denoted in the proof of \cite[Theorem 3.1]{LY201801}.
{\bf Proof}
Recall that $$\overline{X}_t:=-\log\inf\limits_{x\in[a_0,b_0]} \bfP^{x}(\forall_{s\in [0,t]} \beta W_{s}+a\leq B_{s} \leq \beta W_{s}+b,~\beta W_{t}+a'\leq B_{t} \leq \beta W_{t}+b'|W).$$
We first consider the case that
\begin{eqnarray}\label{B2.3.3}
a_0\leq a'< b'\leq b_0,~~~\min\{a_0-a,b-b_0\}>\max\{a'-a_0,b_0-b'\}.
\end{eqnarray}
We still choose constants $a'', b''$ satisfying $a'<a''<b''<b'$ and choose constants $\delta$ such that
$$0<2\beta\delta < \min\big\{b''-a'',~\min\{a_0-a,b-b_0\}-\max\{a''-a_0,b_0-b''\} \big\},$$
After the discussion (3.4)-(3.11) in \cite{LY201801}, we replace the coming $``(\log n)^q"$ in the original proof by $2n.$ Hence the (3.12) in \cite{LY201801} will be rewritten as for any $\varsigma>0$ and large enough $n$ we have
\begin{eqnarray} \label{2.3.13}\bfP\big(\overline{X}_t\geq 2n\big)&\leq& \sum_{i=1}^{+\infty}\bfP\Big(\sum_{k=1}^{i}\frac{C_1}{\rho_{_{k,\delta}}}\geq n,N=i \Big)\no
\\&\leq& \sum_{i=1}^{\lfloor \varsigma n\rfloor-1}\bfP\Big(\sum_{k=1}^{i}\frac{C_1}{\rho_{_{k,\delta}}}\geq n \Big)+\sum_{i=\lfloor \varsigma n\rfloor}^{+\infty}\bfP(N=i)\no
\\&\leq& \varsigma n\bfP\Bigg(\sum_{k=1}^{\lfloor\varsigma n\rfloor}\frac{1}{\rho_{_{k,\delta}}}\geq \frac{n}{C_1}\Bigg)+
\bfP(\tau_{_{\lfloor\varsigma n\rfloor,\delta}}< t),
\end{eqnarray}
where $\rho_{_{k,\delta}}:=\tau_{_{k,\delta}}-\tau_{_{k-1,\delta}}, k\in\bfN^+,$ $\tau_{_{0,\delta}}:=0, \tau_{_{n+1,\delta}}:=\inf\{s>\tau_{_{n,\delta}}: |W_{s}-W_{\tau_{_{n,\delta}}}|=\delta\},~n\in\bfN.$
We show \eqref{bmt+++} by reduction to absurdity. If we assume that there exists $c>0$ such that $\bfE(e^{c\overline{X}_t})=+\infty.$ It is well-known that the  probability density function $p(t)$ of $\tau_{1,\delta}$ is
$$p(t)=\frac{2\delta}{\sqrt{2\pi t^3}}\sum\limits_{n=-\infty}^{+\infty}(4n+1)e^{-\frac{(4n+1)^2\delta^2}{2t}}.$$
Then it is not hard to see for any $c_8>0,$ we have $\bfE(e^{c_8/\tau_{_{1,\delta}}})<+\infty.$ Hence according to Cram\'{e}r theorem (\cite[Page 27]{DZ1998}) we can find a $\varsigma>0$ small enough such that $\bfP\Bigg(\sum_{k=1}^{\lfloor\varsigma n\rfloor}\frac{1}{\rho_{_{k,\delta}}}\geq \frac{n}{C_1}\Bigg)\leq e^{-4cn}$ for large enough $n.$ On the other hand, \cite[Page 30]{IM1974} tells that there exists a constant $c_9>0$  such that $\bfE(e^{c_{_9}\tau_{_{1,\delta}}})<+\infty$, which implies that $\bfE(e^{\lambda\tau_{_{1,\delta}}})\ra 0$ as $\lambda\ra -\infty$ and we can also apply the Cram\'{e}r theorem to estimate $\bfP(\tau_{_{\lfloor\varsigma n\rfloor,\delta}}< t)$. Note that the $t$ is given and not depend on $n,$ hence for the fixed $\varsigma,$ we can find large enough $n$ such that $\bfP(\tau_{_{\lfloor\varsigma n\rfloor,\delta}}< t)\leq e^{-4cn}.$ At last, by the analysis around (3.13) in \cite{LY201801} we can remove the restriction (4.1). In conclusion we get $\bfP\big(\overline{X}_t\geq 2n\big)\leq 2e^{-4cn}$ for large enough $n,$ which contradicts the hypothesis $\bfE(e^{c\overline{X}_t})=+\infty$.

 \ack
We would like to thank Bastien Mallein
for his useful discussions and comments on \cite{LY201801}, which provides a basis for this subject.
This work is supported by the Fundamental Research Funds for the Central Universities (NO.2232021D-30) and the National Natural Science Foundation of China (NO.11971062).

\end{document}